\documentclass[12pt]{article}
\usepackage{amssymb, amsmath, amsfonts, tikz}
\usepackage[indent,margin=1cm]{caption}
\usepackage[colorlinks,linkcolor=blue,anchorcolor=blue,citecolor=blue]{hyperref}
\newtheorem{theo}{Theorem}[section]

\newtheorem{conj}[theo]{Conjecture}

\newtheorem{coro}[theo]{Corollary}

\def\pf{\noindent{\it Proof.} }
\def\qed{\nopagebreak\hfill{\rule{4pt}{7pt}}
\medbreak}

\numberwithin{equation}{section}

\allowdisplaybreaks

\begin{document}

\begin{center}
{\large \bf Stanley's conjectures on the Stern poset}
\end{center}
\begin{center}
Arthur L.B. Yang\\[6pt]

Center for Combinatorics, LPMC\\
Nankai University, Tianjin 300071, P. R. China\\[6pt]

{\tt yang@nankai.edu.cn}
\end{center}

\noindent\textbf{Abstract.} The Stern poset $\mathcal{S}$ is a graded infinite poset naturally associated to Stern's triangle, which was defined by Stanley analogously to Pascal's triangle. Let $P_n$ denote the interval of $\mathcal{S}$ from the unique element of row $0$  of Stern's triangle to the $n$-th element of row $r$ for sufficiently large $r$. For $n\geq 1$ let
\begin{align*}
L_n(q)&=2\cdot\left(\sum_{k=1}^{2^n-1}A_{P_k}(q)\right)+A_{P_{2^n}}(q),
\end{align*}
where $A_{P}(q)$ represents the corresponding $P$-Eulerian polynomial. For any $n\geq 1$ Stanley conjectured that $L_n(q)$ has only real zeros and $L_{4n+1}(q)$ is divisible by $L_{2n}(q)$. In this paper we obtain a simple recurrence relation satisfied by $L_n(q)$ and affirmatively solve Stanley's conjectures. We also establish the asymptotic normality of the coefficients of $L_n(q)$.

\noindent \emph{AMS Classification 2020:} 05A15, 11B37, 26C10


\section{Introduction}

Recently, Stanley \cite{stanley2020acow} introduced a class of polynomials $b_{n}(q)$  by defining
$b_1(q)=1$ and
\begin{align}
b_{2n}(q)&=b_n(q),\label{bPol-1}\\[5pt]
b_{4n+1}(q)&=qb_{2n}(q)+b_{2n+1}(q),\label{bPol-2}\\[5pt]
b_{4n+3}(q)&=b_{2n+1}(q)+qb_{2n+2}(q).\label{bPol-3}
\end{align}
Let
\begin{align}\label{LPol-eq}
L_n(q)&=2\cdot\left(\sum_{k=1}^{2^n-1}b_k(q)\right)+b_{2^n}(q).
\end{align}
The main objective of this paper is to prove Stanley's conjectures on the real-rootedness and divisibility concerning $L_n(q)$.

Let us first review some background. Stanley's conjectures considered here arose in the study of Stern's triangle $S$, which is an array of numbers similar to Pascal's triangle. We follow Stanley \cite{stanley2020aam} to give a description of Stern's triangle. We number the rows of Stern's triangle by consecutive natural numbers beginning with $0$. Row $0$ consists of a single $1$, row $1$ consists of three $1$'s, and for $r\geq 2$ row $r$ is obtained from row $r-1$ by inserting between any two consecutive elements their sums and then placing a $1$ at the beginning and end.
It is clear that row $n$ consists of $2^{r+1}-1$ terms. We number the elements of row $r$ from $0$ to $2^{r+1}-2$ and use
$\left\langle {\genfrac{}{}{0pt}{}{r}{n}}\right\rangle$ to denote the $(n+1)$-th element of row $r$. Thus, we have the following recurrence relation
\begin{align}\label{rec-sterntriangle}
  \left\langle {\genfrac{}{}{0pt}{}{r}{2n+1}}\right\rangle = \left\langle {\genfrac{}{}{0pt}{}{r-1}{n}}\right\rangle, \qquad \left\langle {\genfrac{}{}{0pt}{}{r}{2n}}\right\rangle = \left\langle {\genfrac{}{}{0pt}{}{r-1}{n-1}}\right\rangle+\left\langle {\genfrac{}{}{0pt}{}{r-1}{n}}\right\rangle,
\end{align}
where we set $\left\langle {\genfrac{}{}{0pt}{}{r}{n}}\right\rangle=0$ for $n<0$ or $n>2^{r+1}-2$ for convenience.
Stanley \cite{stanley2020aam} showed that for any $m\geq 1$ the summation
$$\sum_n\left\langle {\genfrac{}{}{0pt}{}{r}{n}}\right\rangle^m$$
 obeys a homogeneous linear recurrence  with constant coefficients, and conjectured the least order of a homogeneous linear recurrence.
Speyer \cite{speyer2019} proved that the above sum satisfies such a recurrence of the conjectured minimal order. 
For any $r\geq 1$ Stanley showed that
\begin{align}\label{eq-gf}
\sum_{n\geq 0} \left\langle {\genfrac{}{}{0pt}{}{r}{n}}\right\rangle x^n=\prod_{i=0}^{r-1} (1+x^{2^i}+x^{2\cdot 2^i}).
\end{align}
Letting $r\rightarrow\infty$ in \eqref{eq-gf}, we get
$$\prod_{i\geq 0} (1+x^{2^i}+x^{2\cdot 2^i})=\sum_{n\geq 1} b_{n} x^{n-1},$$
where the sequence $\{b_{n}\}_{n\geq 0}$ with $b_0=0$ is the well-known Stern's diatomic sequence \cite{stern1858}. For more information on Stern's diatomic sequence, see Northshield \cite{northshield2010}. It is known that $\{b_{n}\}_{n\geq 0}$ satisfies the following recurrence relation
\begin{align}
  b_{2n} & =b_{n}, \qquad b_{2n+1}=b_{n}+b_{n+1}.
\end{align}
Comparing the above recurrence relation with \eqref{bPol-1}, \eqref{bPol-2} and \eqref{bPol-3}, we see that $\{b_{n}(q)\}_{n\geq 1}$ is a polynomial analogue of Stern's diatomic sequence.

\begin{figure}[ht]
\center
\begin{tikzpicture}
[place/.style={thick,fill=black!100,circle,inner sep=0pt,minimum size=2mm,draw=black!100}]

\node [place] (v00) at (8,0) {};
\node at (8.6,0) {${{\left\langle {\genfrac{}{}{0pt}{}{0}{0}}\right\rangle}}$};

\node [place] (v10) at (5,2) {};
\node [place] (v11) at (8,2) {};
\node [place] (v12) at (11,2) {};

\node at (5.6,2) {${{\left\langle {\genfrac{}{}{0pt}{}{1}{0}}\right\rangle}}$};
\node at (8.6,2) {${{\left\langle {\genfrac{}{}{0pt}{}{1}{1}}\right\rangle}}$};
\node at (11.6,2) {${{\left\langle {\genfrac{}{}{0pt}{}{1}{2}}\right\rangle}}$};

\draw [thick](v00) -- (v10);
\draw [thick](v00) -- (v11);
\draw [thick](v00) -- (v12);

\node [place] (v20) at (3.5,4) {};
\node [place] (v21) at (5,4) {};
\node [place] (v22) at (6.5,4) {};
\node [place] (v23) at (8,4) {};
\node [place] (v24) at (9.5,4) {};
\node [place] (v25) at (11,4) {};
\node [place] (v26) at (12.5,4) {};

\node at (4.1,4) {${{\left\langle {\genfrac{}{}{0pt}{}{2}{0}}\right\rangle}}$};
\node at (5.6,4) {${{\left\langle {\genfrac{}{}{0pt}{}{2}{1}}\right\rangle}}$};
\node at (7.1,4) {${{\left\langle {\genfrac{}{}{0pt}{}{2}{2}}\right\rangle}}$};
\node at (8.6,4) {${{\left\langle {\genfrac{}{}{0pt}{}{2}{3}}\right\rangle}}$};
\node at (10.1,4) {${{\left\langle {\genfrac{}{}{0pt}{}{2}{4}}\right\rangle}}$};
\node at (11.6,4) {${{\left\langle {\genfrac{}{}{0pt}{}{2}{5}}\right\rangle}}$};
\node at (13.1,4) {${{\left\langle {\genfrac{}{}{0pt}{}{2}{6}}\right\rangle}}$};

\draw [thick](v10) -- (v20);
\draw [thick](v10) -- (v21);
\draw [thick](v10) -- (v22);

\draw [thick](v11) -- (v22);
\draw [thick](v11) -- (v23);
\draw [thick](v11) -- (v24);

\draw [thick](v12) -- (v24);
\draw [thick](v12) -- (v25);
\draw [thick](v12) -- (v26);

\node [place] (v30) at (2.75,5.5) {};
\node [place] (v31) at (3.5,5.5) {};
\node [place] (v32) at (4.25,5.5) {};
\node [place] (v33) at (5,5.5) {};
\node [place] (v34) at (5.75,5.5) {};
\node [place] (v35) at (6.5,5.5) {};
\node [place] (v36) at (7.25,5.5) {};
\node [place] (v37) at (8,5.5) {};
\node [place] (v38) at (8.75,5.5) {};
\node [place] (v39) at (9.5,5.5) {};
\node [place] (v310) at (10.25,5.5) {};
\node [place] (v311) at (11,5.5) {};
\node [place] (v312) at (11.75,5.5) {};
\node [place] (v313) at (12.5,5.5) {};
\node [place] (v314) at (13.25,5.5) {};

\node at (2.75,6) {${{\left\langle {\genfrac{}{}{0pt}{}{3}{0}}\right\rangle}}$};
\node at (3.5,6) {${{\left\langle {\genfrac{}{}{0pt}{}{3}{1}}\right\rangle}}$};
\node at (4.25,6) {${{\left\langle {\genfrac{}{}{0pt}{}{3}{2}}\right\rangle}}$};
\node at (5,6) {${{\left\langle {\genfrac{}{}{0pt}{}{3}{3}}\right\rangle}}$};
\node at (5.75,6) {${{\left\langle {\genfrac{}{}{0pt}{}{3}{4}}\right\rangle}}$};
\node at (6.5,6) {${{\left\langle {\genfrac{}{}{0pt}{}{3}{5}}\right\rangle}}$};
\node at (7.25,6) {${{\left\langle {\genfrac{}{}{0pt}{}{3}{6}}\right\rangle}}$};
\node at (8,6) {${{\left\langle {\genfrac{}{}{0pt}{}{3}{7}}\right\rangle}}$};
\node at (8.75,6) {${{\left\langle {\genfrac{}{}{0pt}{}{3}{8}}\right\rangle}}$};
\node at (9.5,6) {${{\left\langle {\genfrac{}{}{0pt}{}{3}{9}}\right\rangle}}$};
\node at (10.25,6) {${{\left\langle {\genfrac{}{}{0pt}{}{3}{10}}\right\rangle}}$};
\node at (11,6) {${{\left\langle {\genfrac{}{}{0pt}{}{3}{11}}\right\rangle}}$};
\node at (11.75,6) {${{\left\langle {\genfrac{}{}{0pt}{}{3}{12}}\right\rangle}}$};
\node at (12.5,6) {${{\left\langle {\genfrac{}{}{0pt}{}{3}{13}}\right\rangle}}$};
\node at (13.25,6) {${{\left\langle {\genfrac{}{}{0pt}{}{3}{14}}\right\rangle}}$};


\draw [thick](v20) -- (v30);
\draw [thick](v20) -- (v31);
\draw [thick](v20) -- (v32);

\draw [thick](v21) -- (v32);
\draw [thick](v21) -- (v33);
\draw [thick](v21) -- (v34);

\draw [thick](v22) -- (v34);
\draw [thick](v22) -- (v35);
\draw [thick](v22) -- (v36);

\draw [thick](v23) -- (v36);
\draw [thick](v23) -- (v37);
\draw [thick](v23) -- (v38);

\draw [thick](v24) -- (v38);
\draw [thick](v24) -- (v39);
\draw [thick](v24) -- (v310);

\draw [thick](v25) -- (v310);
\draw [thick](v25) -- (v311);
\draw [thick](v25) -- (v312);

\draw [thick](v26) -- (v312);
\draw [thick](v26) -- (v313);
\draw [thick](v26) -- (v314);

\node at (8,7.5) {$\bullet$};
\node at (8,6.5) {$\bullet$};
\node at (8,7) {$\bullet$};

%
\end{tikzpicture}

\caption{The Stern poset $\mathcal{S}$}\label{fig-sternposet}
\end{figure}
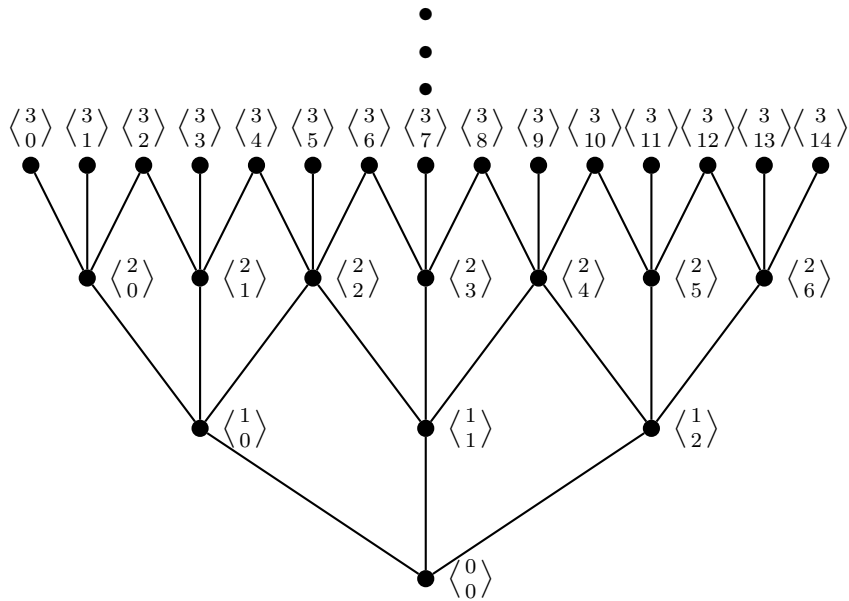

Stanley \cite{stanley2020acow} showed that the polynomials $b_{n}(q)$ also arise as $P$-Eulerian polynomials of certain posets $P$ naturally associated to  Stern's triangle $S$.
To take $S$ as a poset, we will consider $\left\langle {\genfrac{}{}{0pt}{}{r}{n}}\right\rangle$ as a symbol instead of a number.
 According to \eqref{rec-sterntriangle}, we may impose a partial order $\preceq_{\mathcal{S}}$ on $S$ by letting
\begin{align*}
\left\langle {\genfrac{}{}{0pt}{}{r-1}{n}}\right\rangle \preceq_{\mathcal{S}} \left\langle {\genfrac{}{}{0pt}{}{r}{2n+1}}\right\rangle,\qquad
\left\langle {\genfrac{}{}{0pt}{}{r-1}{n}}\right\rangle \preceq_{\mathcal{S}} \left\langle {\genfrac{}{}{0pt}{}{r}{2n}}\right\rangle,\qquad
\left\langle {\genfrac{}{}{0pt}{}{r-1}{n}}\right\rangle \preceq_{\mathcal{S}} \left\langle {\genfrac{}{}{0pt}{}{r}{2n+2}}\right\rangle
\end{align*}
for $0\leq n\leq 2^{r}-2$
and then taking the transitive closure. Following Stanley \cite{stanley2020acow} we call $(S,\, \preceq_{\mathcal{S}})$ the Stern poset, denoted by $\mathcal{S}$. See Figure \ref{fig-sternposet} for the first four levels of the Stern poset.

Fixing a positive integer $n$, suppose that
$$2^{k} \leq n-1< 2^{k+1}$$
for some $k\geq 0$. A little thought shows that for $r\geq k+1$ the interval $\left[\left\langle {\genfrac{}{}{0pt}{}{0}{0}}\right\rangle,\, \left\langle {\genfrac{}{}{0pt}{}{r}{n-1}}\right\rangle\right]$ of $\mathcal{S}$ is the ordinal sum of the chain
$\left[\left\langle {\genfrac{}{}{0pt}{}{0}{0}}\right\rangle,\, \left\langle {\genfrac{}{}{0pt}{}{r-k-1}{0}}\right\rangle\right]$ and the interval $\left[\left\langle {\genfrac{}{}{0pt}{}{r-k-1}{0}}\right\rangle,\, \left\langle {\genfrac{}{}{0pt}{}{r}{n-1}}\right\rangle\right]$, the latter being isomorphic to the interval $\left[\left\langle {\genfrac{}{}{0pt}{}{0}{0}}\right\rangle,\, \left\langle {\genfrac{}{}{0pt}{}{k+1}{n-1}}\right\rangle\right]$ of $\mathcal{S}$ for any $r\geq k+1$. Thus for sufficiently large $r$ we may associate to the $n$-th element of row $r$ in Stern's triangle the poset $\left[\left\langle {\genfrac{}{}{0pt}{}{0}{0}}\right\rangle,\, \left\langle {\genfrac{}{}{0pt}{}{k+1}{n-1}}\right\rangle\right]$, denoted by $P_n$. Stanley \cite{stanley2020acow} obtained the following result.

\begin{theo}[\cite{stanley2020acow}]
For any $n\geq 1$ let $b_n(q)$ be defined as in \eqref{bPol-1}, \eqref{bPol-2} and \eqref{bPol-3}. Then
$b_n(q)$ is equal to the $P_n$-Eulerian polynomial, namely
$$b_n(q)=\sum_{\sigma\in\mathcal{L}(P_n)}q^{\mathrm{des}\,(\sigma)},$$
where $\mathcal{L}(P_n)$ denotes the set of linear extensions of the poset $P_n$, provided that $P_n$ is naturally labeled.
\end{theo}

Therefore, the polynomials $L_n(q)$ defined in \eqref{LPol-eq} can be considered as
Eulerian row sums of Stern's triangle. Stanley \cite{stanley2020acow} proposed the following interesting conjectures.

\begin{conj}[\cite{stanley2020acow}]\label{conj-rz}
For any $n\geq 1$ the polynomial $L_n(q)$ has only real zeros.
\end{conj}

\begin{conj}[\cite{stanley2020acow}]\label{conj-di}
For any $n\geq 1$ the polynomial $L_{4n+1}(q)$ is divisible by $L_{2n}(q)$.
\end{conj}

The aim of this paper is to prove the above conjectures.
Motivated by Conjecture \ref{conj-rz}, we also study the asymptotically normality of $L_n(q)$ since the coefficients of many real-rooted polynomials  follow asymptotically normal limit laws, see \cite{bender1973,harper1967,hwang2020,wangzhangzhu2017} for instance. Recall that, for
a polynomial sequence  $\{F_n(q)\}_{n\geq 0}$ given by
\begin{align}\label{def-f}
F_n(q)=\sum_{k=0}^{n}a_{n,k}q^k
\end{align}
with $a_{n,k}$ all nonnegative, the coefficients $a_{n,k}$ are said to be asymptotically normal by a central limit theorem if
\begin{align}
\lim_{n\rightarrow \infty}\sup_{x\in\mathbb{R}}\left|\sum_{k\leq \mu_n+x\sigma_n}p(n,k)-\frac{1}{\sqrt{2\pi}}\int_{-\infty}^{x}e^{-t^2/2}dt\right|=0,
\end{align}
where
$$
p(n,k)=\frac{a_{n,k}}{\sum_{j=0}^{n}a_{n,j}}
$$
and
$\mu_n,\sigma_n^2$ are respectively the mean and variance of $a_{n,k}$.
We say that $a_{n,k}$ are asymptotically normal by a local limit theorem on $\mathbb{R}$ if
\begin{align}
\lim_{n\rightarrow \infty}\sup_{x\in\mathbb{R}}\left|\sigma_np(n,\lfloor\mu_n+x\sigma_n\rfloor)-\frac{1}{\sqrt{2\pi}}e^{-x^2/2}\right|=0.
\end{align}

The rest of this paper is organized as follows. In Section \ref{sect-rz} we will first derive a recurrence relation satisfied by
$L_n(q)$ and then give a proof of Conjecture \ref{conj-rz}. In Section \ref{sect-di} we will prove Conjecture \ref{conj-di} based on  the recurrence relations satisfied by $L_{2n}(q)$, $L_{2n+1}(q)$ and $L_{4n+1}(q)$. In Section \ref{sect-normal} we will show that the coefficients of
$L_n(q)$ are asymptotically normal by central and local limit theorems.

\section{Real zeros} \label{sect-rz}

The aim of this section is to prove Conjecture \ref{conj-rz}.
Our proof is based on the following recurrence relation satisfied by $L_n(q)$.

\begin{theo}\label{key-rec-thm}
For any $n\geq 2$ we have
\begin{align}\label{key-rec}
L_{n+1}(q)=3L_n(q)+2(q-1)L_{n-1}(q).
\end{align}
\end{theo}

\pf  Substituting \eqref{LPol-eq} into \eqref{key-rec} we obtain that
\begin{align*}
2\cdot\left(\sum_{k=1}^{2^{n+1}-1}b_k(q)\right)+b_{2^{n+1}}(q)=&3\cdot\left(2\sum_{k=1}^{2^n-1}b_k(q)+b_{2^n}(q)\right)\\[5pt]
&+2(q-1)\cdot\left(2\cdot\left(\sum_{k=1}^{2^{n-1}-1}b_k(q)\right)+b_{2^{n-1}}(q)\right).
\end{align*}
By using \eqref{bPol-1} and \eqref{bPol-3}, we find that
\begin{align*}
b_{2^n}(q)=b_2(q)=b_1(q)=1 \mbox{  and } b_3(q)=1+q.
\end{align*}
Thus it suffices to show that
\begin{align*}
2\cdot\left(\sum_{k=4}^{2^{n+1}-1}b_k(q)\right)+2(1+q)+5
=&3\cdot\left(2\sum_{k=2}^{2^n-1}b_k(q)+3\right)\\[5pt]
&+2(q-1)\cdot\left(2\cdot\left(\sum_{k=1}^{2^{n-1}-1}b_k(q)\right)+1\right),
\end{align*}
or equivalently
\begin{align}\label{mid-id}
\sum_{k=4}^{2^{n+1}-1}b_k(q)=3\cdot\left(\sum_{k=2}^{2^n-1}b_k(q)\right)+2(q-1)\cdot\left(\sum_{k=1}^{2^{n-1}-1}b_k(q)\right).
\end{align}
By using \eqref{bPol-1}, \eqref{bPol-2} and \eqref{bPol-3}, it is easy to verify that
\begin{align}\label{key-id}
b_{4m}(q)+b_{4m+1}(q)+b_{4m+2}(q)+b_{4m+3}(q)=(1+q)b_{2m}(q)+3b_{2m+1}(q)+qb_{2m+2}(q).
\end{align}
Substituting \eqref{key-id} into the left hand side of \eqref{mid-id}, we get that
for $n \geq 2$
\begin{align*}
\sum_{k=4}^{2^{n+1}-1}b_k(q)&=\sum_{m=1}^{2^{n-1}-1}\left(b_{4m}(q)+b_{4m+1}(q)+b_{4m+2}(q)+b_{4m+3}(q)\right)\\[5pt]
&=\sum_{m=1}^{2^{n-1}-1}\left((1+q)b_{2m}(q)+3b_{2m+1}(q)+qb_{2m+2}(q)\right)
\end{align*}
By rewriting $(1+q)b_{2m}(q)$ as $3b_{2m}(q)+(q-2)b_{2m}(q)$, we find that
\begin{align*}
\sum_{k=4}^{2^{n+1}-1}b_k(q)&=\sum_{m=1}^{2^{n-1}-1}\left(3b_{2m}(q)+3b_{2m+1}(q)+(q-2)b_{2m}(q)+qb_{2m+2}(q)\right)\\[5pt]
&=3\cdot\left(\sum_{m=1}^{2^{n-1}-1}\left(b_{2m}(q)+b_{2m+1}(q)\right)\right)+\sum_{m=1}^{2^{n-1}-1}\left((q-2)b_{2m}(q)+qb_{2m+2}(q)\right)
\end{align*}
By virtue of \eqref{bPol-1}, we further obtain that
\begin{align*}
\sum_{k=4}^{2^{n+1}-1}b_k(q)=&3\cdot\left(\sum_{k=2}^{2^n-1}b_k(q)\right)+\sum_{m=1}^{2^{n-1}-1}\left((q-2)b_{m}(q)+qb_{m+1}(q)\right)\\[5pt]
=&3\cdot\left(\sum_{k=2}^{2^n-1}b_k(q)\right)+q\cdot\left(\sum_{m=1}^{2^{n-1}-1}b_{m}(q)\right)\\[5pt]
&+q\cdot \left(\sum_{m=1}^{2^{n-1}-1}b_{m+1}(q)\right)-2\cdot \left(\sum_{m=1}^{2^{n-1}-1}b_{m}(q)\right).
\end{align*}
Note that, for the third summation on the right hand side, we have
\begin{align*}
\sum_{m=1}^{2^{n-1}-1}b_{m+1}(q)&=b_{2^{n-1}}(q)+\sum_{m=1}^{2^{n-1}-2}b_{m+1}(q)\\[5pt]
&=b_{1}(q)+\sum_{m=2}^{2^{n-1}-1}b_{m}(q)\\[5pt]
&=\sum_{m=1}^{2^{n-1}-1}b_{m}(q)
\end{align*}
where the second equality holds since $b_{2^{n-1}}(q)=b_1(q)=1$ by \eqref{bPol-1}.
Thus
\begin{align*}
\sum_{k=4}^{2^{n+1}-1}b_k(q)=3\cdot\left(\sum_{k=2}^{2^n-1}b_k(q)\right)+2(q-1)\cdot\left(\sum_{k=1}^{2^{n-1}-1}b_k(q)\right),
\end{align*}
as desired in \eqref{mid-id}. This completes the proof.
\qed

In view of \eqref{key-rec} the polynomial sequence $\{L_n(q)\}_{n\geq 1}$ lies in the framework of polynomial sequence of type $(0,1)$ studied by Gross, Mansour, Tucker and Wang \cite{gross2016}, who obtained the following result.

\begin{theo}[{\cite[Theorem 2.6]{gross2016}}]\label{gross-thm}
Let $\{M_n(q)\}_{n\geq 0}$ be the polynomial sequence defined by the recursion
\begin{align*}
  M_{n+1}(q) & =aM_{n}(q)+(bq+c)M_{n-1}(q)
\end{align*}
with initial values $M_0(q)=1$ and $M_1(q)=t(q-r)$, where $a,b,t>0, c,r\in \mathbb{R}$, and $r\neq -c/b$.
Then for any $n\geq 1$ the polynomial $M_n(q)$ has only real zeros.
\end{theo}

The main result of this section is as follows, which gives an affirmative answer to Conjecture \ref{conj-rz}.

\begin{theo}\label{main-thm1}
For any $n\geq 1$ the polynomial $L_n(q)$ has only real zeros.
\end{theo}

\pf By \eqref{LPol-eq} it is straightforward to compute that
$L_1(q)=3$ and $L_{2}(q)=2q+7$. For any $n\geq 0$ let
\begin{align*}
 \tilde{L}_n(q) & =\frac{{L}_{n+1}(q)}{3}.
\end{align*}
It is clear that
$$\tilde{L}_{n+1}(q)=3\tilde{L}_n(q)+2(q-1)\tilde{L}_{n-1}(q)$$
with $\tilde{L}_1(q)=1$ and $\tilde{L}_{2}(q)=\frac{2}{3}q+\frac{7}{3}$.
It is easy to verify that $\{\tilde{L}_n(q)\}_{n\geq 0}$ satisfies the conditions of Theorem \ref{gross-thm}.
Thus for each $n\geq 0$ the polynomial $\tilde{L}_n(q)$ has only real zeros, so does ${L}_{n+1}(q)$.
\qed

\noindent\textit{Remark.} Based on the recurrence relation \eqref{key-rec}, Theorem \ref{main-thm1} can also be proved by using a result due to Wang and Yeh, see \cite[Theorem 1]{wangyeh2005} or \cite[Theorem 1.1]{liuwang2007}.

Theorem \ref{key-rec-thm} also has some other consequences. For instance, we can determine the ordinary generating function of $\{L_{n}(q)\}_{n\geq 1}$.
Let
$$
\Phi(x)=\sum_{n\geq 1} L_{n}(q) x^n.
$$
By \eqref{key-rec} it is easy to show that
\begin{align*}
\Phi(x)&=\frac{x(2(q-1)x+3)}{1-3x-2(q-1)x^2}=\frac{x(2(q-1)x+3)}{1-3x-2(q-1)x^2}=\frac{1}{1-x(2(q-1)x+3)}-1.
\end{align*}
We can also determine the shifted Hankel determinants of $\{L_{n}(q)\}_{n\geq 1}$.
Let
$$\mathcal{H}_n^{(k)}=(L_{i+j+k-1}(q))_{1\leq i,j\leq n},\quad {H}_n^{(k)}=\det \mathcal{H}_n^{(k)}.$$
We obtain the following result.

\begin{coro}
For any $k\geq 0$ and $m\geq 3$ we have
$${H}_1^{(k)}=L_{k+1}(q), \quad {H}_2^{(k)}=(-1)^{k+1}2^{k+2}(q-1)^{k+2}, \quad {H}_m^{(k)}=0.$$
\end{coro}

\pf We only need to compute ${H}_2^{(k)}$ and ${H}_m^{(k)}$ for $m\geq 3$.
It is routine to compute that
\begin{align*}
{H}_2^{(0)}&=\det
\begin{pmatrix}
L_{1}(q) & L_{2}(q)\\[5pt]
L_{2}(q) & L_{3}(q)
\end{pmatrix}=\det
\begin{pmatrix}
3 & 2q+7\\[5pt]
2q+7 & 12q+15
\end{pmatrix}=-2^2(q-1)^2.
\end{align*}
For $k\geq 1$ applying \eqref{key-rec} to ${H}_2^{(k)}$ gives
\begin{align}
{H}_2^{(k)}&=\det
\begin{pmatrix}
L_{k+1}(q) & L_{k+2}(q)\\[5pt]
L_{k+2}(q) & L_{k+3}(q)
\end{pmatrix}\nonumber\\[8pt]
&=\det
\begin{pmatrix}
L_{k+1}(q) & 3L_{k+1}(q)+2(q-1)L_{k}(q) \\[5pt]
L_{k+2}(q) & 3L_{k+2}(q)+2(q-1)L_{k+1}(q)
\end{pmatrix}\nonumber\\[8pt]
&=\det
\begin{pmatrix}
L_{k+1}(q) & 2(q-1)L_{k}(q) \\[5pt]
L_{k+2}(q) & 2(q-1)L_{k+1}(q)
\end{pmatrix}\nonumber\\[8pt]
&=2(q-1)\cdot \det
\begin{pmatrix}
L_{k+1}(q) & L_{k}(q) \\[5pt]
L_{k+2}(q) & L_{k+1}(q)
\end{pmatrix}\nonumber\\[8pt]
&=-2(q-1)\cdot \det
\begin{pmatrix}
L_{k}(q) & L_{k+1}(q) \\[5pt]
L_{k+1}(q) & L_{k+2}(q)
\end{pmatrix}\nonumber\\[8pt]
&=-2(q-1){H}_2^{(k-1)}.\label{rec-hankel}
\end{align}
By iterating \ref{rec-hankel} we obtain
\begin{align*}
{H}_2^{(k)}=(-1)^k2^k(q-1)^k{H}_2^{(0)}=(-1)^{k+1}2^{k+2}(q-1)^{k+2},
\end{align*}
as desired.

We proceed to compute ${H}_m^{(k)}$ for $m\geq 3$. By definition
\begin{align*}
\mathcal{H}_m^{(k)}=\begin{pmatrix}
L_{k+1}(q) & L_{k+2}(q) & L_{k+3}(q) & \cdots & L_{k+m}(q)\\[5pt]
L_{k+2}(q) & L_{k+3}(q) & L_{k+4}(q) & \cdots & L_{k+m+1}(q)\\[5pt]
\vdots & \vdots & \vdots & \cdots & \vdots\\[5pt]
L_{k+m}(q) & L_{k+m+1}(q) & L_{k+m+2}(q) & \cdots & L_{k+2m-1}(q)
\end{pmatrix}.
\end{align*}
By \eqref{key-rec} the third column of $\mathcal{H}_m^{(k)}$ is a linear combination of the first column and the second column, and thus the Hankel determinant ${H}_m^{(k)}$ vanishes. This completes the proof.
\qed

\section{Divisibility} \label{sect-di}

In this section we aim to prove Conjecture \ref{conj-di}. The idea of the proof is to introduce a polynomial sequence $\{M_n(q)\}_{n\geq 1}$ as defined below, and then to show that $\{M_n(q)L_{2n}(q)\}_{n\geq 1}$ and  $\{L_{4n+1}(q)\}_{n\geq 1}$ satisfy the same recurrence relation with the same initial conditions.

We first give recurrence relations which are satisfied by $\{L_{2n}(q)\}_{n\geq 1}$, $\{L_{2n-1}(q)\}_{n\geq 1}$ and $\{L_{4n+1}(q)\}_{n\geq 1}$ respectively.
For notational convenience let
\begin{align}\label{eq-eoj}
E_n(q)&=L_{2n}(q), \qquad O_n(q)=L_{2n-1}(q), \qquad J_n(q)=L_{4n+1}(q)
\end{align}
 for any $n\geq 1$.
The recurrence relations are as follows.

\begin{theo}\label{thm-di}
Let
$$
a=(5+4q), b=-4(q-1)^2, h=2O_3(q). 
$$
Then for any $n\geq 2$ we have
\begin{align}
E_{n+1}(q)&=aE_{n}(q)+bE_{n-1}(q),\label{rec-e}\\[5pt]
O_{n+1}(q)&=aO_{n}(q)+bO_{n-1}(q),\label{rec-o}\\[5pt]
J_{n+1}(q)&=ab^{n-1}h+a^2J_{n}(q)+b^2J_{n-1}(q)+\sum_{k=2}^{n-1}2a^2b^{n-k}J_k(q).\label{rec-j}
\end{align}
\end{theo}

\pf Let us first prove \eqref{rec-e}.
From \eqref{key-rec} it follows that
\begin{align*}
E_{n+1}(q)=&L_{2(n+1)}(q)=3L_{2n+1}(q)+2(q-1)L_{2n}\\[5pt]
=&3(3L_{2n}(q)+2(q-1)L_{2n-1}(q))\\[5pt]
&+2(q-1)(3L_{2n-1}(q)+2(q-1)L_{2(n-1)}(q)).
\end{align*}
By simplification we obtain
\begin{align*}
E_{n+1}(q)&=9L_{2n}(q)+12(q-1)L_{2n-1}(q)+4(q-1)^2L_{2(n-1)}(q),
\end{align*}
which could be rewritten as
\begin{align*}
E_{n+1}(q)=&(5+4q)L_{2n}(q)-4(q-1)^2L_{2(n-1)}(q)\\[5pt]
&+4(1-q)\left(L_{2n}(q)-3L_{2n-1}(q)-2(q-1)L_{2(n-1)}(q)\right).
\end{align*}
Again by \eqref{key-rec} we obtain
\begin{align*}
E_{n+1}(q)&=(5+4q)E_{n}(q)-4(q-1)^2E_{n-1}(q),
\end{align*}
as desired.

Similarly, we compute
\begin{align*}
O_{n+1}(q)&=L_{2(n+1)+1}(q)\\[5pt]
&=3L_{2(n+1)}(q)+2(q-1)L_{2n+1}\\[5pt]
&=3(3L_{2n+1}(q)+2(q-1)L_{2n}(q))+2(q-1)(3L_{2n}(q)+2(q-1)L_{2n-1}(q))\\[5pt]
&=9L_{2n+1}(q)+12(q-1)L_{2n}(q)+4(q-1)^2L_{2n-1}(q)\\[5pt]
&=(5+4q)L_{2n+1}(q)-4(q-1)^2L_{2(n-1)+1}(q)\\[5pt]
&\quad +4(1-q)\left(L_{2n+1}(q)-3L_{2n}(q)-2(q-1)L_{2n-1}(q)\right)\\[5pt]
&=(5+4q)O_{n}(q)-4(q-1)^2O_{n-1}(q),
\end{align*}
where the last equality follows from \eqref{key-rec}. This completes the proof of \eqref{rec-o}.

We proceed to prove \eqref{rec-j}.  Note that
\begin{align*}
J_{n+1}(q)&=L_{4(n+1)+1}(q)=O_{2(n+1)}(q).
\end{align*}
By virtue of \eqref{rec-o}, there holds
\begin{align}
J_{n+1}(q)&=aO_{2n+1}(q)+bO_{2n}(q)\nonumber\\[5pt]
&=a\left(aO_{2n}(q)+bO_{2n-1}(q)\right)+b\left(aO_{2n-1}(q)+bO_{2n-2}(q)\right)\nonumber\\[5pt]
&=a^2J_n(q)+b^2J_{n-1}(q)+2abO_{2n-1}(q).\label{rec-jo}
\end{align}
On the other hand, for any $m\geq 2$ we have
\begin{align*}
O_{2m-1}(q)&=aO_{2m-2}(q)+bO_{2m-3}(q)=aJ_{m-1}(q)+bO_{2m-3}(q).
\end{align*}
Iterating the above process leads to
\begin{align}
O_{2n-1}(q)&=aJ_{n-1}(q)+bO_{2n-3}(q)\nonumber\\[5pt]
&=aJ_{n-1}(q)+b\left( aJ_{n-2}(q)+bO_{2n-5}(q) \right)\nonumber\\[5pt]
&=aJ_{n-1}(q)+b\left( aJ_{n-2}(q)+b\left( aJ_{n-3}(q)+bO_{2n-7}(q)\right) \right)\nonumber\\[5pt]
&=\cdots\nonumber\\[5pt]
&=b^{n-2}O_3(q)+\sum_{k=2}^{n-1}ab^{n-k-1}J_k(q).\label{rec-mid}
\end{align}
Substituting \eqref{rec-mid} into \eqref{rec-jo}, we get
\begin{align*}
J_{n+1}(q)&=a^2J_n(q)+b^2J_{n-1}(q)+2ab\left(b^{n-2}O_3(q)+\sum_{k=2}^{n-1}ab^{n-k-1}J_k(q)\right)\\[5pt]
&=2ab^{n-1}O_3(q)+a^2J_{n}(q)+b^2J_{n-1}(q)+\sum_{k=2}^{n-1}2a^2b^{n-k}J_k(q),
\end{align*}
as desired. The proof is complete.
\qed

We continue to prove Conjecture \ref{conj-di}. Now define the sequence $\{M_n(q)\}_{n\geq 1}$ by
letting $M_1(q)=\frac{J_1(q)}{E_1(q)}$ and $M_2(q)=\frac{J_2(q)}{E_2(q)}$
and
\begin{align}\label{rec-m}
M_{n+1}(q)&=aM_{n}(q)+bM_{n-1}(q),
\end{align}
where $a=(5+4q)$ and $b=-4(q-1)^2$ as in Theorem \ref{thm-di}.
The main result of this section is as follows, which gives an affirmative answer to Conjecture \ref{conj-di}.

\begin{theo}
For any $n\geq 1$  let $E_n(q), J_n(q)$ be given by \eqref{eq-eoj}, let $M_n(q)$ be given by \eqref{rec-m}, and let $H_n(q)$ be defined by
\begin{align}
H_n(q)=E_n(q)M_{n}(q).\label{def-h}
\end{align}
Then $H_n(q)=J_n(q)$, and hence $L_{4n+1}(q)$ is divisible by $L_{2n}(q)$.
\end{theo}

\pf
Since $H_1(q)=J_1(q)$ and $H_2(q)=J_2(q)$, it suffices to show that the sequence $\{H_n(q)\}_{n\geq 1}$ satisfies the same recurrence relation as
 $\{J_n(q)\}_{n\geq 1}$, namely
\begin{align}
H_{n+1}(q)&=ab^{n-1}h+a^2H_{n}(q)+b^2H_{n-1}(q)+\sum_{k=2}^{n-1}2a^2b^{n-k}H_k(q),\label{rec-h}
\end{align}
where
$$
a=(5+4q), b=-4(q-1)^2, h=2O_3(q). 
$$

By \eqref{rec-e} and \eqref{rec-m} we have
\begin{align*}
  H_{n+1}(q) & =E_{n+1}(q)M_{n+1}(q)\\[5pt]
  & =\left(aE_{n}(q)+bE_{n-1}(q)\right)\left(aM_{n}(q)+bM_{n-1}(q)\right)\\[5pt]
  & =a^2H_n(q)+b^2H_{n-1}(q)+ab\left(E_{n}(q)M_{n-1}(q)+E_{n-1}(q)M_{n}(q)\right).
\end{align*}
For notational convenience, let
\begin{align}\label{eq-polt}
T_n(q)=E_{n}(q)M_{n-1}(q)+E_{n-1}(q)M_{n}(q).
\end{align}
Then by \eqref{rec-e} and \eqref{rec-m}
\begin{align*}
T_n(q)&=E_{n}(q)M_{n-1}(q)+E_{n-1}(q)M_{n}(q)\\
&=\left(aE_{n-1}(q)+bE_{n-2}(q)\right)M_{n-1}(q)+E_{n-1}(q)\left(aM_{n-1}(q)+bM_{n-2}(q)\right)\\
&=2aH_{n-1}(q)+bT_{n-1}(q),
\end{align*}
where the last equality is due to \eqref{def-h} and \eqref{eq-polt}.

Iterating the above process, we get
\begin{align*}
  H_{n+1}(q) & =a^2H_n(q)+b^2H_{n-1}(q)+abT_n(q)\\[5pt]
  & =a^2H_n(q)+b^2H_{n-1}(q)+ab\left(2aH_{n-1}(q)+bT_{n-1}(q)\right)\\[5pt]
  & =a^2H_n(q)+b^2H_{n-1}(q)+2a^2bH_{n-1}(q)+ab^2T_{n-1}(q)\\[5pt]
  & =a^2H_n(q)+b^2H_{n-1}(q)+2a^2bH_{n-1}(q)+ab^2\left(2aH_{n-2}(q)+bT_{n-2}(q)\right)\\[5pt]
  & =a^2H_n(q)+b^2H_{n-1}(q)+2a^2bH_{n-1}(q)+2a^2b^2 H_{n-2}(q)+ab^3T_{n-2}(q)\\[5pt]
  &=\cdots\\[5pt]
  & =a^2H_n(q)+b^2H_{n-1}(q)+2a^2bH_{n-1}(q)+2a^2b^2 H_{n-2}(q)+\cdots\\[5pt]
  &\quad + 2a^2b^{n-2} H_{2}(q)+ab^{n-1}T_{2}(q)\\
  & =a^2H_n(q)+b^2H_{n-1}(q)+\sum_{k=2}^{n-1}2a^2b^{n-k}H_k+ab^{n-1}T_{2}(q).
\end{align*}
It is straightforward to verify that
$$T_2(q)=2O_3(q).$$
Thus $\{H_{n}(q)\}_{n\geq 1}$ and $\{J_{n}(q)\}_{n\geq 1}$ satisfy the same recurrence relation with the same initial conditions.
This completes the proof.
\qed

\section{Asymptotic normality}\label{sect-normal}

In this section we will show that the coefficients of
$L_n(q)$ are asymptotically normal by central and local limit theorems.
For a polynomial sequence  $\{F_n(q)\}_{n\geq 0}$ as in \eqref{def-f},
it is not hard to show that if the coefficients $a_{n,k}$ are asymptotically normal by a local limit theorem on $\mathbb{R}$ then these numbers are also asymptotically normal by a central limit theorem. Thus, we will focus on proving that the coefficients of
$L_n(q)$ are asymptotically normal by a local limit theorem.

The main tool we use here is the following criterion, which goes back to Harper \cite{harper1967}.

\begin{theo}[{\cite[Theorem 2]{bender1973}}]\label{lemm-as}
Suppose that $\{F_n(q)\}_{n\geq 0}$ is a real-rooted polynomial sequence with nonnegative coefficients as in \eqref{def-f}.
Let
\begin{align}
\mu_n&=\frac{F_n'(1)}{F_n(1)},\label{eq-mean}\\[5pt]
\sigma_n^2&=\frac{F_n''(1)}{F_n(1)}+\mu_n-\mu_n^2.\label{eq-var}
\end{align}
If $\sigma_n^2\rightarrow +\infty$ when $n \rightarrow +\infty$, then the coefficients of $F_n(q)$ are asymptotically normal by a local limit theorem with mean $\mu_n$ and variance $\sigma_n^2$.
\end{theo}

In order to compute the mean and variance corresponding to $L_n(q)$, we next derive Binet's formula for these polynomials by solving  \eqref{key-rec}. Note that the characteristic equation of this recurrence relation is
\begin{align*}
x^2-3x-2(q-1)=0
\end{align*}
with roots
$$
r(q)=\frac{3+\sqrt{8q+1}}{2} \quad \mbox{and}\quad s(q)=\frac{3-\sqrt{8q+1}}{2}
$$
where $r(q)+s(q)=3$ and $r(q)s(q)=-2(q-1)$. So the general solution of \eqref{key-rec} is
$$
L_n(q)=C(r(q))^{n-1}+D(s(q))^{n-1},
$$
where the coefficients $C$ and $D$ are to be determined.

The initial conditions $L_1(q)=3$ and $L_2(q)=7+2q=3r(q)+3s(q)-r(q)s(q)$ yield the following system:
\begin{align*}
C+D&=3,\\[5pt]
Cr(q)+Ds(q)&=2q+7.
\end{align*}
Solving this, we get
$$C=\frac{2q+7-3s(q)}{r(q)-s(q)} \quad \mbox{and}\quad  D=-\frac{2q+7-3r(q)}{r(q)-s(q)}.$$
Thus
\begin{align}\label{eq-Binet}
L_n(q)&=\frac{2q+7-3s(q)}{r(q)-s(q)}\cdot (r(q))^{n-1}-\frac{2q+7-3r(q)}{r(q)-s(q)}\cdot (s(q))^{n-1}\nonumber\\[5pt]
&=\frac{3-s(q)}{r(q)-s(q)}\cdot (r(q))^{n}-\frac{3-r(q)}{r(q)-s(q)}\cdot (s(q))^{n}\nonumber\\[5pt]
&=3\cdot \frac{(r(q))^{n}-(s(q))^{n}}{r(q)-s(q)}-r(q)s(q)\cdot \frac{(r(q))^{n-1}-(s(q))^{n-1}}{r(q)-s(q)}.
\end{align}

The main result of this section is as follows.

\begin{theo}
Let the sequence $\{L_n(q)\}_{n\geq 1}$ be defined as in \eqref{LPol-eq}.
Then the coefficients of $L_n(q)$ are asymptotically normal by a local limit theorem.
\end{theo}

\pf Let us first compute $\mu_n$ and $\sigma_n^2$.  It is easy to show that
\begin{align*}
r(q)&=\frac{3+\sqrt{8q+1}}{2},\\[5pt]
r'(q)&=\frac{2}{\sqrt{8q+1}},\\[5pt]
r''(q)&=-\frac{8}{(8q+1)\sqrt{8q+1}},\\[5pt]
s(q)&=\frac{3-\sqrt{8q+1}}{2},\\[5pt]
s'(q)&=-\frac{2}{\sqrt{8q+1}},\\[5pt]
s''(q)&=\frac{8}{(8q+1)\sqrt{8q+1}}.
\end{align*}
Thus
\begin{align}
r(1)=3,\quad r'(1)=\frac{2}{3},\quad r''(1)=-\frac{8}{27},\\
s(1)=0, \quad s'(1)=-\frac{2}{3}, \quad s''(1)=\frac{8}{27}.
\end{align}
Let
$$F_n(q)=\frac{(r(q))^{n}-(s(q))^{n}}{r(q)-s(q)}.$$
It is straightforward to verify that
\begin{align*}
F'_n(q)=&\frac{n[((r(q))^{n-1}r'(q)-(s(q))^{n-1}s'(q))(r(q)-s(q))]}{(r(q)-s(q))^2}\\[5pt]
&-\frac{[((r(q))^{n}-(s(q))^{n})(r'(q)-s'(q))]}{(r(q)-s(q))^2}
\end{align*}
and
\begin{align*}
F''_n(q)=\frac{G_n'(q)(r(q)-s(q))^2-2G_n(q)(r(q)-s(q))(r'(q)-s'(q))}{(r(q)-s(q))^4}.
\end{align*}
with
\begin{align*}
G_n(q)=&(n-1)((r(q))^{n}r'(q)+(s(q))^{n}s'(q))\\[5pt]
&+((s(q))^{n-1}-n(r(q))^{n-1})r'(q)s(q)
+((r(q))^{n-1}-n(s(q))^{n-1})s'(q)r(q)
\end{align*}
and
\begin{align*}
G_n'(q)=&(n-1)((r(q))^{n}r''(q)+n(r(q))^{n-1}(r'(q))^2+(s(q))^{n}s''(q)+n(s(q))^{n-1}(s'(q))^2)\\[5pt]
&+((s(q))^{n-1}-n(r(q))^{n-1})(r''(q)s(q)+r'(q)s'(q))\\[5pt]
&+((n-1)(s(q))^{n-2}s'(q)-n(n-1)(r(q))^{n-2}r'(q))r'(q)s(q)\\[5pt]
&+((r(q))^{n-1}-n(s(q))^{n-1})(s''(q)r(q)+s'(q)r'(q))\\[5pt]
&+((n-1)(r(q))^{n-2}r'(q)-n(n-1)(s(q))^{n-2}s'(q))s'(q)r(q).
\end{align*}
From the above formulas we deduce that
\begin{align*}
F_n(1)&=3^{n-1},\quad F'_n(1)=2(n-2)3^{n-3},\quad G_n(1)=2(n-2)3^{n-1},\\[5pt]
G_n'(1)&=4(n^2-3n+4)3^{n-3}, \quad F_n''(1)=4(n^2-7n+12)3^{n-5}.
\end{align*}

By Binet's formula for $L_n(q)$, we obtain
\begin{align*}
L_n(q)=3F_n(q)-r(q)s(q)F_{n-1}(q)
\end{align*}
and hence
\begin{align*}
L_n'(q)&=3F_n'(q)-r'(q)s(q)F_{n-1}(q)-r(q)s'(q)F_{n-1}(q)-r(q)s(q)F'_{n-1}(q)\\[5pt]
L_n''(q)&=3F_n''(q)-r''(q)s(q)F_{n-1}(q)-r'(q)s'(q)F_{n-1}(q)-r'(q)s(q)F'_{n-1}(q)\\[5pt]
&-r'(q)s'(q)F_{n-1}(q)-r(q)s''(q)F_{n-1}(q)-r(q)s'(q)F'_{n-1}(q)\\[5pt]
&-r'(q)s(q)F'_{n-1}(q)-r(q)s'(q)F'_{n-1}(q)-r(q)s(q)F''_{n-1}(q).
\end{align*}
So
\begin{align*}
L_n(1)&=3^n,\quad L_n'(1)=2(n-1)\cdot 3^{n-2}, \quad L''_n(1)=4(n-2)(n-3)\cdot 3^{n-4}.
\end{align*}

From \eqref{eq-mean} it follows that
\begin{align*}
  \mu_n&=\frac{L_n'(1)}{L_n(1)}=\frac{2(n-1)}{9}
\end{align*}
and
\begin{align*}
\sigma_n^2&=\frac{L_n''(1)}{L_n(1)}+\mu_n-\mu_n^2\\
&=\frac{4(n-2)(n-3)}{81}+\frac{2(n-1)}{9}-\frac{4(n-1)^2}{81}\\
&=\frac{2(n-1)(2n+7)}{81},
\end{align*}
which tends to $+\infty$ as $n$ approaches infinity.
Then combining Theorem \ref{main-thm1} and Theorem \ref{lemm-as}, we complete the proof.
 \qed

\noindent \textbf{Acknowledgments.}
This work is supported in part by the Fundamental Research Funds for the Central Universities and the National Science Foundation of China (Nos. 11522110, 11971249).

\end{document}